%% file: CPR_explicit__arXiv.tex
\documentclass{scrartcl}

\usepackage[utf8]{luainputenc}
\usepackage[USenglish]{babel}
\usepackage{csquotes}

\usepackage[a4paper, top=2.5cm, bottom=2.5cm, left=2.5cm, right=2.5cm]{geometry}
\usepackage[plainpages=false,pdfpagelabels,hidelinks,unicode]{hyperref}

\usepackage[%
  backend=biber,
  style=numeric-comp,
  giveninits=true, uniquename=init, 
  natbib=true,
  url=true,
  doi=true,
  isbn=false,
  backref=false,
  maxnames=99,
  ]{biblatex}
\usepackage{amsmath}
\allowdisplaybreaks
\usepackage{amssymb}
\usepackage{commath}
\usepackage{mathtools}
\usepackage{bbm}
\usepackage{nicefrac}

\usepackage{amsthm}
\theoremstyle{plain}
  \newtheorem{theorem}{Theorem}
  \newtheorem{lemma}[theorem]{Lemma}
  
\theoremstyle{definition}
  
  \newtheorem{remark}[theorem]{Remark}
  \newtheorem{example}[theorem]{Example}
\numberwithin{theorem}{section}

\usepackage{siunitx}

\usepackage{color}
\usepackage{graphicx}
\usepackage[small]{caption}
\usepackage{subcaption}

\begingroup\expandafter\expandafter\expandafter\endgroup
\expandafter\ifx\csname pdfsuppresswarningpagegroup\endcsname\relax
\else
\fi

\usepackage{booktabs}
\usepackage{rotating}
\usepackage{multirow}

\usepackage{ifluatex}
\ifluatex
  \usepackage[no-math]{fontspec}
\else
  \usepackage[T1]{fontenc}
\fi
\usepackage{newpxtext,newpxmath}

\usepackage{calc}
\usepackage{xparse}

\renewcommand{\vec}[1]{\underline{#1}}
\NewDocumentCommand{\mat}{mo}{%
  \IfValueTF{#2}{%
    \underline{\underline{#1}}{#2}
  }{%
    \underline{\underline{#1}}\,
  }%
}
\newcommand{\diag}[1]{\operatorname{diag}\left(#1\right)}

\newcommand{\scp}[2]{\left\langle{#1,\, #2}\right\rangle}

\renewcommand{\L}{\mathbf{L}}
\renewcommand{\P}{\mathbb{P}}

\newcommand{\fnum}{f^{\mathrm{num}}}
\newcommand{\vecfnum}{\vec{f}^{\mathrm{num}}}

\renewcommand{\epsilon}{\varepsilon}
\renewcommand{\phi}{\varphi}
\newcommand{\N}{\mathbb{N}}
\newcommand{\R}{\mathbb{R}}

\newsavebox{\DelimiterBox}
\newlength{\DelimiterHeight}
\newlength{\DelimiterDepth}
\newsavebox{\ArgumentBox}
\newlength{\ArgumentHeight}
\newlength{\ArgumentDepth}
\newlength{\ResizedDelimiterHeight}

\newenvironment{keywords}{\par\textbf{Key words.}}{\par}
\newenvironment{AMS}{\par\textbf{AMS subject classification.}}{\par}

\title{Analysis of Artificial Dissipation of Explicit and Implicit Time-Integration Methods}
\author{Philipp \"Offner, Jan Glaubitz, Hendrik Ranocha}
\date{December 6, 2019}

\makeatletter
\hypersetup{pdfauthor={\@author}}
\hypersetup{pdftitle={\@title}}
\makeatother

\begin{document}

\maketitle

\begin{abstract}
Stability is an important aspect of numerical methods for hyperbolic conservation laws and has received much interest. However, continuity in time is often assumed and only semidiscrete stability is studied. Thus, it is interesting to investigate the influence of explicit and implicit time integration methods on the stability of numerical schemes. If an explicit time integration method is applied, spacially stable numerical schemes for hyperbolic conservation laws can result in unstable fully discrete schemes. Focusing on the explicit Euler method (and convex combinations thereof), undesired terms in the energy balance trigger this phenomenon and introduce an erroneous growth of the energy over time. In this work, we study the influence of artificial dissipation and modal filtering in the context of discontinuous spectral element methods to remedy these issues. In particular, lower bounds on the strength of both artificial dissipation and modal filtering operators are given and an adaptive procedure to conserve the (discrete) $\L_2$ norm of the numerical solution in time is derived. This might be beneficial in regions where the solution is smooth and for long time simulations. Moreover, this approach is used to study the connections between explicit and implicit time integration methods and the associated energy production. By adjusting the adaptive procedure, we demonstrate that filtering in explicit time integration methods is able to mimic the dissipative behavior inherent in implicit time integration methods. This contribution leads to a better understanding of existing algorithms and numerical techniques, in particular the application of artificial dissipation as well as modal filtering in the context of numerical methods for hyperbolic conservation laws together with the selection of explicit or implicit time integration methods.
\end{abstract}



\input{1_Introduction}

\input{2_FR_SBP}

\input{3_adaptive_filtering}
\input{4_explicit_implicit}

\input{5_Summary}

\appendix

\input{Appendix}

\section*{Acknowledgments}
Jan Glaubitz was supported by the German Research Foundation (DFG, Deutsche
Forschungsgemeinschaft) under Grant SO 363/15-1.
Philipp \"Offner was supported by SNF project ``Solving advection dominated
problems with high order schemes with polygonal meshes: application to compressible
and incompressible flow problems'' and the UZH Postdoc Grant.
Hendrik Ranocha was supported by the German Research Foundation (DFG,
Deutsche Forschungsgemeinschaft) under Grant SO~363/14-1.

\printbibliography

\end{document}

%% file: 1_Introduction.tex
\section{Introduction}

Stability is one of the main desirable properties for a numerical scheme
to solve hyperbolic conservation laws.
This is due to the fact that at least for linear symmetric systems,
an energy estimate (and the correct number of boundary conditions for initial boundary
value problems) comes along with uniqueness and existence of a solution
\cite{gustafsson1995time}.
In the last years, several approaches have been proposed to construct
entropy stable/conservative schemes like in
\cite{tadmor1987numerical, tadmor2003entropy, ranocha2016summation, ranocha2017shallow, offner2019stability, 
gassner2016well, wintermeyer2015entropy, abgrall2018some, abgrall2018general, chen2017entropy, carpenter2014entropy,
fjordholm2012arbitrarily}
and references therein.
Recently, Abgrall \cite{abgrall2018general}
presented a way to build entropy stable/conservative schemes
using the Residual Distribution (RD) framework.
In \cite{abgrall2018connection},
this idea is extended to Flux Reconstruction (FR) methods.
This idea is fairly general and has been extended and re-interpreted in the discontinuous Galerkin (DG) context
in \cite{abgrall2019reinterpretation}.
However, besides the spacial discretization, the selection of the time
integration method is essential for stability of these methods.

First of all, one has to choose between explicit or
implicit methods to march forward in time.
Implicit methods have favorable stability properties and, in particular, allow larger time steps.
For instance, by using implicit time integration methods build on Summation-By-Parts (SBP) operators in
time\footnote{These schemes can be interpreted
as implicit Runge-Kutta (RK) methods \cite{boom2015high, ranocha2019some}.} \cite{nordstrom2013summation},
the semidiscrete stability results transfer directly to the fully discrete case
\cite{nikkar2015fully, nordstrom2017roadmap, friedrich2018entropy}.
It should be stressed, however, that implicit methods yield to (typically non-)linear systems to be solved.
Since the time step is also constrained by accuracy requirements, implicit methods
may not be as efficient as explicit ones.

Explicit time integration methods, on the other hand, can be faster and easier to implement, but suffer
from stability issues triggered by additional error terms.
One way to improve the stability properties of numerical schemes
is the usage of artificial dissipation.
This idea dates back to early works of von Neumann and Richtymer \cite{vonneumann1950method}.
Since then, many researchers have contributed to the development and extension of artificial dissipation methods,
including the works
\cite{tadmor1987numerical, tadmor1989convergence, ma1998chebyshev, mattsson2004stable, ranocha2018stability}.

In this work, we investigate the connections between artificial dissipation in explicit time
integration methods and implicit time integration methods without additional limiting from point of stability.
We further extend this investigation to modal filtering.
Modal filtering is strongly connected to artificial dissipation methods in spectral and spectral element methods
\cite{majda1978fourier, glaubitz2018application, kreiss1979stability, gottlieb2001spectral, canuto2006spectral, offner2013spectral} 
and provides an alternative which, 
in some cases\footnote{For instance, when the method is already formulated in a suitable modal basis.},
might be more efficient and easier to implement. 
In particular, we demonstrate that it is possible to mimic the dissipation (and thus stability) inherent in implicit
time integration methods for explicit time integration methods when modal filtering with a suitable filter strength is
incorporated.
This result directly carries over to explicit time integration methods with suitable artificial dissipation terms.
Thus, we are able to present an approach to obtain stable fully discrete schemes using explicit time
integration.
Such discretizations combine the favorable stability properties of implicit time integration methods with the
efficiency gain of explicit time integration methods.
Finally, we would like to mention that recently a relaxation Runge-Kutta approach 
has been proposed to construct fully discrete explicit energy (entropy) conservative/stable schemes
in \cite{ketcheson2019relaxation, ranocha2019relaxation}. Their approach
has some similarities to our consideration but instead of working with modal filters or artificial viscosity 
to destroy the energy production in time, they change the final update step in the RK method to guarantee 
that the discrete energy equality is fulfilled.

For sake of simplicity, the explicit Euler method is considered.
Yet, at least for non-linear problems, the same stability issues arise for strong stability preserving (SSP)
RK schemes, since they can be written as convex combination of explicit Euler steps \cite{gottlieb2011strong}.
In the appendix, we show how our investigation carries over to general Runge-Kutta methods.
Recent relevant articles concerned with the strong stability of explicit Runge-Kutta
methods are, e.g., \cite{sun2017stability, sun2018strong, lozano2018entropy,
lozano2018entropyImplicit, ranocha2018strong, ranocha2018L2stability}.

The rest of this work is organized as follows:
In section \ref{sec:FR}, we start by briefly revisiting the FR method in its formulation using SBP operators.
This method yields a stable semidiscretisation and thus serves as a representative of a stable scheme.
Yet, the examinations are rather general and valid for other spacial discretizations as well.
In section \ref{sec:artificial-dissipation}, we investigate the mechanism which triggers
stability issues when semidiscretisation (even stable ones) are evolved in time by explicit time marching.
Further, we investigate the stabilizing effect of artificial dissipation terms and modal filtering.
In principle, similar investigations are well-known.
Performing this analysis in a vector matrix-vector representation including suitable discrete inner products, however,
we are able derive new (strict) bounds on the artificial viscosity strength and filter strength for stability to carry
over in time.
Building up on this strategy, adaptive filtering strategies can be derived which yield methods with neither not enough
nor too much dissipation.
This might be beneficial in smooth regions for long time simulations.
Section \ref{sec:exp_imp} explores the connection between implicit time integration and modal filtering in explicit
time integration.
We end this work by a summary in section \ref{Summary}.
In the appendix, we demonstrate how the investigation for the filter strength of section
\ref{sec:artificial-dissipation} can be extended to general Runge-Kutta methods.

%% file: 2_FR_SBP.tex
\section{Flux Reconstruction using Summation-By-Parts Operators}
\label{sec:FR}

In this section, we provide a brief description of FR
methods using SBP operators, which will serve as a representative for spacially stable methods in the later
investigations.
Yet, it should be stressed that our analysis is also valid for other space discretisation, such as
like DG or finite volume (FV) schemes.
Further, let us consider a one-dimensional scalar conservation law
\begin{equation}
\label{eq:scalar-CL}
  \partial_t u + \partial_x f(u) = 0
\end{equation}
on $\Omega = [a,b]$, equipped with adequate boundary and initial conditions.
For sake of simplicity, in this work, periodic boundary conditions will be assumed.

The domain is partitioned into smaller subdomains, also called elements, which are mapped diffeomorphically
onto a reference element, typically $\Omega_{\mathrm{ref}} = [-1,1]$.
All calculations are conducted within this reference element then.
There, the solution $u$ is approximated by a polynomial $U$ of degree $ \leq N$.
Let $\{ \zeta_i \}_{i=0}^N$ be a set of interpolation points in $[-1,1]$.
Then, $U \in \P_N$ can be represented w.r.t. to a nodal Lagrange basis, resulting in a vector of coefficients $\vec{u}$
given by $\vec{u}_i = U(\zeta_i)$ for $i = 0, \dots, N$.

Note that the representation of $U$ w.r.t. other bases is possible as well.
In the setting described in \cite{ranocha2017extended}, the solution is represented
as an element of a finite dimensional Hilbert space of functions on the volume.
W.r.t. a chosen basis, the scalar product approximating the $L_2$ scalar
product is represented by a matrix $\mat{M}$ and the derivative (divergence) by
$\mat{D}$. Additionally, functions on the boundary (consisting of two points in
this one dimensional case) are elements of another finite dimensional Hilbert
space with appropriate basis. The restriction of functions on the volume to the
boundary is represented by a (rectangular) matrix $\mat{R}$ and integration
w.r.t. the outer normal by $\mat{B} = \diag{-1,1}$. Finally, the operators have to
fulfil the SBP property
\begin{equation}
\label{eq:SBP}
  \mat{M} \mat{D} + \mat{D}[^T] \mat{M}
  = \mat{R}[^T] \mat{B} \mat{R}
\end{equation}
as a compatibility condition in order to mimic integration by parts
\begin{equation}
  \vec{u}^T \mat{M} \mat{D} \vec{v} + \vec{u}^T \mat{D}[^T] \mat{M} \vec{v}
  \approx
  \int_\Omega u \, \partial_x v + \int_\Omega \partial_x u \, v
  = u \, v \big|_{\partial \Omega}
  \approx
  \vec{u}^T \mat{R}[^T] \mat{B} \mat{R} \vec{v}.
\end{equation}
Additional ingredients of FR methods are numerical fluxes (Riemann solvers) $\fnum$,
computing a single valued flux on the boundary using values from both neighbouring
elements, and a correction step which can be formulated
as a Simultaneous-Approximation-Term (SAT) from finite difference (FD) methods \cite{ranocha2016summation}. An overview and translation rules linking the notation
used in this article and in DG or FD methods can be found in \cite{ranocha2018generalised}.

In the following, either nodal Gauß-Legendre and Lobatto-Legendre or modal Legendre bases will be used.
Multiplication of functions on the volume will be conducted pointwise for
nodal bases or exactly, followed by an $L_2$ projection, for modal bases. The
resulting multiplication operators are written with two underlines, e.g.
$\mat{u}$ represents multiplication with the polynomial given by $\vec{u}$.

\begin{example}
 We give two examples for the discretisation. The linear advection
 with constant velocity is given by
\begin{equation}
\label{eq:lin-adv}
  \partial_t u + \partial_x u = 0.
\end{equation}
The semidiscretisation using the canonical choice for the correction procedure can be
written as
\begin{equation}
\label{eq:lin-adv-semidisc}
  \partial_t \vec{u}
  = - \mat{D} \vec{u}
    - \mat{M}[^{-1}] \mat{R}[^T] \mat{B} \left(
        \vecfnum - \mat{R} \vec{u}
      \right)
\end{equation}
in the reference domain. 
The second example, which we also consider later is the
nonlinear Burgers' equation
\begin{equation}
\label{eq:Burgers}
  \partial_t u + \partial_x \frac{u^2}{2} = 0.
\end{equation}
This is more difficult, since discontinuities may develop in finite time.
A discretisation is given by a skew symmetric form
\begin{equation}
\label{eq:Burgers-semidisc}
  \partial_t \vec{u}
  =
  - \frac{1}{3} \mat{D} \mat{u} \vec{u}
  - \frac{1}{3} \mat{u}{^*} \mat{D} \vec{u}
  + \mat{M}[^{-1}] \mat{R}[^T] \mat{B} \left(
      \vecfnum
      - \frac{1}{3} \mat{R} \mat{u} \vec{u}
      - \frac{1}{6} \left( \mat{R} \vec{u} \right)^2
    \right).
\end{equation}
Using \eqref{eq:lin-adv-semidisc} or \eqref{eq:Burgers-semidisc} results in spacially stable schemes if an appropriate
numerical flux is applied, see \cite{ranocha2016summation}.
\end{example}

%% file: 3_adaptive_filtering.tex
\section{Artificial Dissipation and Modal Filtering }
\label{sec:artificial-dissipation}

In this section, we investigate the stabilising effect of artificial dissipation operators and modal filtering.
We note that both techniques share a strong connection.
Using a first order operator splitting in time, artificial dissipation operators can be interpreted as
exponential modal filter, see
\cite{offner2013spectral,glaubitz2018application,ranocha2018stability,glaubitz2019smooth}.


In artificial dissipation methods, a small (super) diffusive term of even order is added to the conservation
law \eqref{eq:scalar-CL}.
This yields
\begin{equation}
 \label{eq:scalar-CL-RHS}
  \partial_t u(t,x) + \partial_x f \left( u(t,x) \right)
  = (-1)^{s+1} \epsilon \left( \partial_x a(x) \partial_x \right)^{s} u(t,x),
 \end{equation}
where $s \in \N$ is the \emph{order}, $\epsilon \geq 0$ the \emph{strength}, and
$a \colon \R \to \R$ is a suitable function.
The term $\left( \partial_x a(x) \partial_x \right)^{s}$ describes the $s$-fold
application of the linear operator
$f(x) \mapsto \partial_x \left( a(x) \partial_x f(x) \right)$. 
Henceforth, the dependence on $t$ and $x$ will be implied but not written explicitly in all cases.


\subsection{Discrete Formulation}

In order to get a working numerical scheme, a time discretisation has to be
introduced.
For simplicity, we start by considering an explicit Euler method.
Yet, once stability is ensured for the simple explicit Euler method, this result carries over to the whole class of
explicit SSP-RK methods under appropriate time step restriction.
This is, for instance, described in the monograph \cite{gottlieb2011strong} and references
cited therein.

In the standard element, one time step $\Delta t$ by the explicit Euler method is
given by
\begin{equation}
\label{eq:exp-Euler}
  \vec{u}
  \mapsto \vec{u}_+ := \vec{u} + \Delta t \, \partial_t \vec{u}.
\end{equation}
Using an SBP-FR semidiscretisation to compute the time derivative
$\partial_t \vec{u}$ in \eqref{eq:exp-Euler} without
artificial viscosity, the norm after one Euler step is given by
\begin{equation}
\label{eq:exp-Euler_2}
\begin{aligned}
  \norm{ \vec{u}_+ }_M^2
  & = \vec{u}_+^T \mat{M} \vec{u}_+ \\
  & = \vec{u}^T \mat{M} \vec{u}
    + 2 \Delta t \, \vec{u}^T \mat{M} \partial_t \vec{u}
    + (\Delta t)^2 \partial_t \vec{u}^T \mat{M} \partial_t \vec{u} \\
  & = \norm{ \vec{u} }_M^2
    + 2 \Delta t \, \scp{ \vec{u} }{ \partial_t \vec{u} }_M
    + (\Delta t)^2 \norm{ \partial_t \vec{u} }_M^2.
\end{aligned}
\end{equation}
Here, the second term on the right hand side,
$2 \Delta t \, \scp{ \vec{u} }{ \partial_t \vec{u} }_M$,
yields only boundary terms that can be controlled by the numerical flux.
However, the last term, $(\Delta t)^2 \norm{ \partial_t \vec{u} }_M^2 \geq 0$,
is non-negative and might therefore increase the norm.
It is this term, which is responsible for (spacially stable) methods to still become unstable in time.
In the following subsections, we investigate two approaches to remedy this source of instability.

\subsection{Application of Artificial Viscosity}
\label{subsec:Application}

We now derive a lower bound on the strength $\epsilon$ for artificial dissipation to carry spacial stability of a
method over to time.
Assuming a fixed function $a$ and order $s$, the strength $\epsilon$ can be
estimated in the following way.
Denoting the time derivative obtained by the underlying SBP-FR method without artificial dissipation by
$\partial_t \vec{u}$ and the matrix representation of the discretised artificial dissipation by
\begin{equation}\label{eq:RHS-smart}
  \mat{A}[^s]
  :=
  \left(\mat{M}[^{-1}] \mat{D}[^T] \mat{M} \mat{a} \mat{D} \right)^s,
\end{equation}
yields
\begin{equation}
\label{eq:av-derivative}
  \partial_t \vec{u}^\epsilon
  =
  \partial_t \vec{u}
  - \epsilon \left( \mat{M}[^{-1}] \mat{D}[^T] \mat{M} \mat{a} \mat{D} \right)^s \vec{u}.
\end{equation}
Note that other discretisations of the artificial viscosity term in \eqref{eq:scalar-CL-RHS} are possible but not
recommended.
Yet, it has been proved in \cite{ranocha2018stability} that the discretisation \eqref{eq:RHS-smart} is compatible with
SBP operators and results in dissipation of the $\L^2$-entropy.
Thus, after one time step by the explicit Euler method \emph{with artificial dissipation}, the norm is given by
\begin{equation}
\label{eq:av-euler-condition}
\begin{aligned}
  \norm{ \vec{u}^\epsilon_+ }_M^2
  =&
  \norm{ \vec{u} }_M^2
  + 2 \Delta t \scp{ \vec{u} }{ \partial_t \vec{u}^\epsilon }_M
  + (\Delta t)^2 \norm{ \partial_t \vec{u}^\epsilon }_M^2
  \\=&
  \norm{ \vec{u} }_M^2
  + 2 \Delta t \scp{ \vec{u} }{ \partial_t \vec{u} }_M
  - 2 \epsilon \Delta t \scp{\vec{u}}{\mat{A}[^s] \vec{u}}_M
  + (\Delta t)^2 \norm{ \partial_t \vec{u}^\epsilon }_M^2.
\end{aligned}
\end{equation}
Again, $\scp{ \vec{u} }{ \partial_t \vec{u} }_M$ can be estimated in terms of
boundary values and numerical fluxes and is negative (non-positive) for a spacially stable discretisation of 
\eqref{eq:scalar-CL}.
Hence, for the method to be stable in time, the two last terms need to cancel out.
In this case,
\begin{equation}
  \norm{ \vec{u}^\epsilon_+ }_M^2
  =
  \norm{ \vec{u} }_M^2
  + 2 \Delta t \scp{ \vec{u} }{ \partial_t \vec{u} }_M
\end{equation}
would follow and the fully discrete scheme will be conservative as well as stable in space \emph{and time}.
Using \eqref{eq:av-euler-condition}, the condition of the last to terms to cancel out can be rewritten as
\begin{equation}
\begin{aligned}
 0
 =&
 - 2 \epsilon \scp{\vec{u}}{\mat{A}[^s] \vec{u}}_M
 + \Delta t \norm{ \partial_t \vec{u}^\epsilon }_M^2
 \\
 =&
 - 2 \epsilon \scp{\vec{u}}{\mat{A}[^s] \vec{u}}_M
 + \Delta t \left( \norm{ \partial_t \vec{u} }_M^2
 - 2 \epsilon \scp{\partial_t \vec{u}}{\mat{A}[^s] \vec{u}}_M
 + \epsilon^2 \norm{ \mat{A}[^s] \vec{u} }_M^2 \right),
\end{aligned}
\end{equation}
which again is equivalent to
\begin{equation}
\label{eq:quadratic}
  \epsilon^2 \underbrace{\left( \Delta t \norm{ \mat{A}[^s] \vec{u} }_M^2 \right)}_{=: X}
  + \epsilon \underbrace{\left( - 2 \scp{\vec{u}}{\mat{A}[^s] \vec{u}}_M
    - 2 \Delta t \scp{\partial_t \vec{u}}{\mat{A}[^s] \vec{u}}_M \right)}_{=: Y}
  + \underbrace{\left( \Delta t \norm{ \partial_t \vec{u} }_M^2 \right)}_{=: Z}
  = 0.
\end{equation}
The (possibly complex) roots of the equation \eqref{eq:quadratic} for $X \neq 0$ are given by
\begin{equation}
  \epsilon_{1/2} = \frac{1}{2 X} \left( -Y \pm \sqrt{ Y^2 -4 X Z } \right).
\end{equation}
Hence, for a sufficiently small time step $\Delta t$ and if the solution is not constant, the discriminant $Y^2 - 4 X
Z$ is non-negative and there is at least one real solution $\epsilon$.
Additionally, both $-Y$ and $X Z$ are positive for
sufficiently small $\Delta t$, since the artificial dissipation operator
$\mat{A}$ is positive semi-definite, i.e.
\begin{equation}
\label{eq:av-condition-Delta-t}
  Y^2 - 4 X Z > 0
  ,\quad
  -Y          > 0
  ,\qquad \text{ if } \Delta t \text{ is small enough and } \mat{A}[^s] \vec{u} \neq 0.
\end{equation}
Thus,
\begin{equation}
  \epsilon_1
  \geq
  \epsilon_2
  =
  \frac{1}{2 X} \left( -Y - \sqrt{ Y^2 -4 X Z } \right)
  \geq
  \frac{1}{2 X} \left( -Y + \sqrt{ Y^2 } \right)
  =
  0,
\end{equation}
and the roots of the quadratic equation \eqref{eq:quadratic} are non-negative.
These results are summed up in the following
\begin{lemma}
\label{lem:discrete}
  If a conservative and stable SBP-FR method for a scalar conservation law
  $\partial_t u + \partial_x f(u) = 0$ 
  is augmented with the artificial dissipation
  $- \epsilon \left( \mat{M}[^{-1}] \mat{D}[^T] \mat{M} \mat{a} \mat{D} \right)^s \vec{u}$ 
  on the right hand side, the fully discrete scheme using an explicit Euler
  method as time discretisation is both conservative and stable if
  \begin{itemize}
    \item
    a nodal Gauß-Legendre / Lobatto-Legendre or a modal Legendre basis is used,

    \item
    $\scp{ \vec{u} }{ \mat{A}[^s] \vec{u} } > 0$, which will be fulfilled for
    the choice of $a$ described below if the solution $\vec{u}$ is not constant,

    \item
    the time step $\Delta t$ is small enough such that \eqref{eq:av-condition-Delta-t}
    is fulfilled,

    \item
    and the strength $\epsilon > 0$ is big enough.
  \end{itemize}

  If the other conditions are fulfilled, $\epsilon$ has to satisfy
  \begin{equation}
  \label{eq:av-epsilon-adaptive}
    \epsilon
    \geq \epsilon_2
    = \frac{1}{2 X} \left( -Y - \sqrt{ Y^2 -4 X Z } \right),
  \end{equation}
  where $X, Y,$ and $Z$ from equation \eqref{eq:quadratic} are used.
\end{lemma}
In our implementation, the strength of dissipation is chosen as the second
 (smaller) root $\epsilon_2$ and results in methods with highly desired
stability properties, as we presented in numerical tests at the end of this section.

\begin{remark}
It remains an interesting, yet unanswered, question how to interpret
the existence of an additional solution $\epsilon_1$. Since this solution yields
a larger strength, the resulting methods show higher dissipation, which might
be undesired in elements without discontinuities or for long time simulations
\cite{oeffner2018error, offner2019error}.
\end{remark}

Note that the CFL condition and therefore the time step in an explicit time integration method depends on the
parameters of the viscosity term.
If no care is taken, artificial dissipation operators will decreases the allowable time step size; 
see \cite{guermond2011entropy,klockner2011viscous,glaubitz2019smooth} and references therein.
Additionally, equation \eqref{eq:av-condition-Delta-t} limits the maximal time
step and can be used as an adaptive strategy
to control this quantity.  This could be also used for an adaptive control strategy and will be
considered in future investigations.
Here, a simple limiting strategy is used for the numerical
experiments.
If the time step is not small enough and equation
\eqref{eq:av-condition-Delta-t} is not fulfilled, the strength $\epsilon$
computed from \eqref{eq:av-epsilon-adaptive} might be negative.
In this case, to avoid instabilities, $\epsilon$ is set to zero, i.e. no artificial viscosity
is used in the corresponding elements.
This phenomenon is strongly connected with stability requirements of the artificial dissipation operator, which have
been discussed above.
Considering a time step by the explicit Euler method for the equation $\partial_t \vec{u} = - \epsilon \mat{A}[^s]
\vec{u}$,
the norm after one time step satisfies
\begin{equation}
\begin{aligned}
  \norm{\vec{u}_+}_M^2
  &=
  \norm{\vec{u}}_M^2
  - 2 \, \epsilon \, \Delta t \scp{ \vec{u} }{ \mat{A}[^s] \vec{u} }_M
  + \epsilon^2 (\Delta t)^2 \norm{\mat{A}[^s] \vec{u}}_M^2
  \\
  &\leq
  \norm{\vec{u}}_M^2
  - 2 \, \epsilon \, \Delta t \scp{ \vec{u} }{ \mat{A}[^s] \vec{u} }_M
  + \epsilon^2 (\Delta t)^2 \norm{\mat{A}}^s \norm{\vec{u}}_M.
\end{aligned}
\end{equation}
Thus, in order to guarantee $\norm{\vec{u}_+}_M^2 \leq \norm{\vec{u}}_M^2$,
for $\mat{A}[^s] \vec{u} \neq 0$, $\Delta t$ has to be limited by
\begin{equation}
  \Delta t
  \leq
  \frac{ 2 \scp{ \vec{u} }{ \mat{A}[^s] \vec{u} }_M }
       { \epsilon \norm{\mat{A}[^s] \vec{u}}_M \norm{\mat{A}}^s \norm{\vec{u}}_M }
  \leq
  \frac{ 2 \norm{\vec{u}}_M \norm{\mat{A}[^s] \vec{u}}_M }
       { \epsilon \norm{\mat{A}[^s] \vec{u}}_M \norm{\mat{A}}^s \norm{\vec{u}}_M }
  =
  \frac{ 2 }{ \epsilon \norm{\mat{A}}^s }.
\end{equation}
Since $\mat{A}$ is a second-order derivative operator, this yields a restriction on the
time step of order $\mathcal{O}\left( \Delta x^{2s} / \epsilon \right)$.
However, it should be noted that $\epsilon$ is computed using the given value of $\Delta t$
and is typically small.

Since Theorem 1 in \cite{ranocha2018stability} requires $a \big|_{[-1,1]} \geq 0$ to be a
polynomial fulfilling $a(\pm 1) = 0$, a simple choice is $a(x) = 1 - x^2$.
By this choice, the continuous artificial dissipation operators is related
to the eigenvalue equation of Legendre polynomials and resulting implications
and connections with modal filtering are  presented in \cite{ranocha2018stability}.

\subsection{Usage of Modal Filters}
\label{subsec:filtering}

In this subsection, we investigate stability of the explicit Euler method combined with modal filtering, which is
strongly related to artificial dissipation 
\cite{majda1978fourier,kreiss1979stability,gottlieb2001spectral,canuto2006spectral,glaubitz2018application,
ranocha2018stability,glaubitz2019smooth}.  
In certain cases, for instance when the method is already formulated w.r.t. a suitable modal basis, 
modal filtering can be more efficient and easier to implement than artificial dissipation.
Further, no additional time step restrictions are introduced.
For modal filtering, an operator splitting approach is applied together with
an explicit Euler method. The update reads
 \begin{equation}
\label{eq:filtered-euler}
  \vec{u}
  \mapsto \tilde{\vec{u}}_+ := \vec{u} + \Delta t \, \partial_t \vec{u}
  \mapsto \vec{u}_+ := \mat{F} \tilde{\vec{u}}_+,
\end{equation}
where \eqref{eq:exp-Euler} holds for $\tilde{\vec{u}}_+$
instead of $\vec{u}_+$ and $\mat{F}$ is the modal filter.
If the filter $\mat{F}$ reduces the norm of $\tilde{\vec{u}}_+$ by the
amount of the additional term $(\Delta t)^2 \left(\partial_t \vec{u} \right)^T
\mat{M} \partial_t \vec{u}$, the fully discrete scheme allows the same estimate
as the semidiscrete one.
Therefore, similar to artificial dissipation, the modal filter has to eliminate the additional positive term.
This idea is summarised in
\begin{lemma}
\label{lem:operator-splitting}
  Rendering a conservative and stable semidiscretisation of the scalar
conservation
  law  \eqref{eq:scalar-CL}
  fully discrete by using an explicit Euler step with modal filtering
  \eqref{eq:filtered-euler} yields a conservative and stable scheme if
  \begin{equation}
  \label{eq:filtered-euler-condition}
  \begin{aligned}
    \norm{\mat{F} \tilde{\vec{u}}_+}_M^2
    &=
    \norm{\vec{u}}_M^2 + 2 \Delta t \scp{\vec{u}}{\partial_t \vec{u}}_M
    \\
    &\leq
    \norm{ \tilde{\vec{u}}_+ }_M^2
    =
    \norm{\vec{u}}_M^2 + 2 \Delta t \scp{\vec{u}}{\partial_t \vec{u}}_M
    + (\Delta t)^2 \norm{\partial_t \vec{u}}_M^2.
  \end{aligned}
  \end{equation}
  This condition can be fulfilled (per element) if
  \begin{itemize}
    \item
    the rate of change $\partial_t \vec{u}$ is zero or

    \item
    the intermediate value $\tilde{\vec{u}}_+$ is not constant and the time step
    $\Delta t$ is small enough.
  \end{itemize}
\end{lemma}
In order to fulfil condition \eqref{eq:filtered-euler-condition} of Lemma
\ref{lem:operator-splitting}, the filter strength $\epsilon$ (with time step
$\Delta t$ included) has to be adapted.
Using a modal Legendre basis, the (exact) modal filter
$\mat{F}$ can be written as
\begin{equation}\label{eq:modal_filter}
  \mat{F} = \diag{ \exp \left[ - \epsilon \, \lambda_n^s \, \Delta t
\right]_{n=0}^p },
\end{equation}
where $\lambda_n = n (n+1) \geq 0$ as it is derived in \cite{ranocha2018stability}.
For stability, the selection of the  free parameter $\epsilon$ is essential.
Similar to subsection \ref{subsec:Application}, we now derive a lower bound on the filter strength that ensures
stability.
Representing the polynomial given by $\tilde{\vec{u}}$ in a modal Legendre basis, i.e. as a
linear combination of Legendre polynomials $\phi_n$, condition
\eqref{eq:filtered-euler-condition} translates to
\begin{equation}
  \sum_{n=0}^p \exp[- 2 \epsilon \, \lambda_n^s \, \Delta t] \,
\tilde{u}_{+,n}^2 \, \norm{\phi_n}^2
  = RHS,
\end{equation}
where the right-hand side $\norm{u}_M^2 + 2 \Delta t
\scp{\vec{u}}{\partial_t\vec{u}}$
is abbreviated as $RHS$.
Using the well-known inequality
\begin{equation}\label{eq:approximation}
  \exp[x] \geq 1 + x, \quad x \in \R,
\end{equation}
as a first order approximation, $\epsilon$ can be estimated by
\begin{equation}
\begin{aligned}
  & \sum_{n=0}^p (1 - 2 \epsilon \, \lambda_n^s \, \Delta t) \,
\tilde{u}_{+,n}^2 \, \norm{\phi_n}^2
    \leq RHS
  \\ \iff
  & \left( \sum_{n=0}^p \tilde{u}_{+,n}^2 \, \norm{\phi_n}^2 - RHS \right)
    \left( \sum_{n=0}^p 2 \lambda_n^s \, \tilde{u}_{+,n}^2 \, \norm{\phi_n}^2
\right)^{-1}
    \leq \epsilon \, \Delta t
\end{aligned}
\end{equation}
for $\sum_{n=0}^p 2 \lambda_n^s \, \tilde{u}_{+,n}^2 \, \norm{\phi_n}^2 > 0$.
Note that we have $\sum_{n=0}^p 2 \lambda_n^s \, \tilde{u}_{+,n}^2 \, \norm{\phi_n}^2 > 0$ if and only if $\tilde{u}_+$
is not identically zero.
Inserting
\begin{equation}
\begin{aligned}
  \sum_{n=0}^p \tilde{u}_{+,n}^2 \, \norm{\phi_n}^2
  &=
  \norm{\vec{u}}_M^2 + 2 \Delta t \scp{\vec{u}}{\partial_t \vec{u}}_M
  + (\Delta t)^2 \norm{\partial_t \vec{u}}_M^2
  \\&=
  RHS + (\Delta t)^2 \norm{\partial_t \vec{u}}_M^2,
\end{aligned}
\end{equation}
this yields
\begin{lemma}
\label{lem:estimate-epsilon}
  A necessary condition for the filter strength according to Lemma
  \ref{lem:operator-splitting} is
  \begin{equation}
  \label{eq:estimate-epsilon}
    \epsilon
    \geq
    \Delta t \norm{\partial_t \vec{u}}_M^2
    \left( \sum_{n=0}^p 2 \lambda_n^s \, \tilde{u}_{+,n}^2 \, \norm{\phi_n}^2
\right)^{-1}.
  \end{equation}
\end{lemma}
\begin{remark}\label{re:Projection_remark}
By applying estimation \eqref{eq:estimate-epsilon} in our numerical scheme, an adaptive strategy
can be applied.
Note that other approximations than \eqref{eq:approximation} could be used.
The same is true if, instead of the explicit Euler method, an (explicit) SSP time integration method is applied, since
such methods can be written as a convex combinations of steps by the explicit Euler method \cite{gottlieb2011strong},
and the filter is applied after each Euler step.
An extension to some Deferred Correction (DEC) methods can  also be done, since one
can write some of theses methods likewise as convex combinations of Euler steps  \cite{liu2008strong, abgrall2017high}.
However, since the triangle inequality is invoked for the resulting estimates,
an undesired additional decrease of the norm may result.
Therefore, the adaptive modal filtering should be applied only after a full time step and not for every stage.
This was for instance demonstrated in \cite{ranocha2018stability}.
Further, this renders  the computation more efficient.
Nevertheless, an extension of this approach to classical Runge-Kutta methods can be done,
yielding to some further conditions which we present in the appendix \ref{sec:appendix}.
\end{remark}

\begin{figure}[!htb]
\centering
  \includegraphics[width=\textwidth]{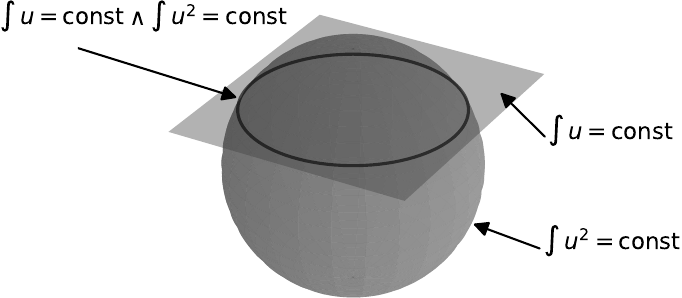}
  \caption{Visualisation of the requirements for projections such as filtering.}
  \label{fig:projection}
\end{figure}

Finally, it should be stressed that adaptive modal filtering can be interpreted as a
special case of projection, enforcing the constraint on the squared norm (a quadratic form) and not violating
conservation, i.e. a constraint on the integral of the solution (a linear form).
This is visualised in Figure \ref{fig:projection}.
However, there are various possibilities to conduct this projection.
As noted in section IV.4 of \cite{hairer2006geometric}, projection methods can
be useful, but can also destroy good properties.
Therefore, they have to be investigated thoroughly.

\begin{remark}[Connection to Relaxation RK approach]
In our analysis, we describe the production of energy in equation \eqref{eq:exp-Euler_2}.
Here, artificial viscosity or modal filters are applied to remove the additional energy.
In
 \cite{ketcheson2019relaxation, ranocha2019relaxation}, the idea is instead to manipulate the time step such
 that the energy remains constant.
 This can be interpreted as a projection in the direction of the step
 update which conserves the energy and all linear invariants. 
 In this article, we use different kinds of projections, e.g. ones given by modal filters,
 which also preserve important linear invariants such as the total mass.
\end{remark}

\subsection{Numerical Simulations}\label{subsec:numerics}

We close this section with a numerical demonstration of the above results and derived adaptive filtering strategies.


\subsubsection*{Comparing Modal Filtering and Projection}

As an example, the linear advection equation with constant coefficients
\begin{equation}\label{eq:transport}
  \partial_t u + \partial_x u = 0
\end{equation}
in $[-1,1]$ with periodic boundary conditions is considered.
For the spacial discretisation, we choose a grid of $N=8$ elements using
polynomials of degree $\leq p = 9$ and an upwind numerical flux.

\begin{figure}[!htb]
  \centering
  \begin{subfigure}[t]{0.49\textwidth}
    \includegraphics[width=\textwidth]{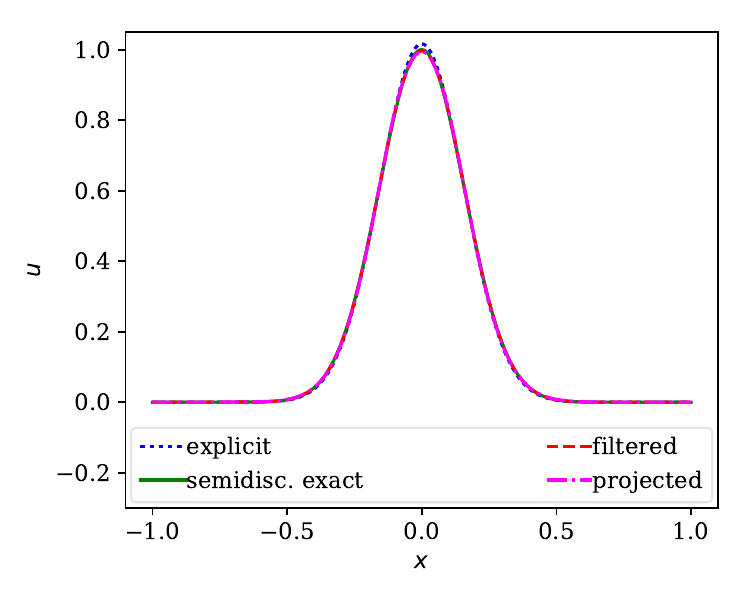}
    \caption{Smooth initial condition \eqref{eq:proj-filt-Euler-exp}.}
    \label{fig:proj-filt-Euler-exp}
  \end{subfigure}%
  ~
  \begin{subfigure}[t]{0.49\textwidth}
    \includegraphics[width=\textwidth]{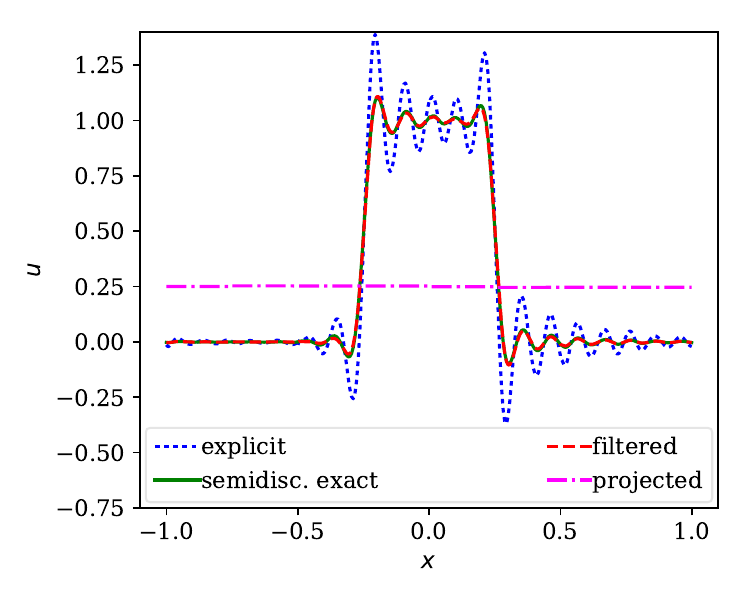}
    \caption{Discontinuous initial condition \eqref{eq:proj-filt-Euler-jump}.}
    \label{fig:proj-filt-Euler-jump}
  \end{subfigure}%
  \caption{Solutions at $t = 4$ computed using \num{20000} time steps of the
           unmodified, filtered, and projected explicit Euler method.}
  \label{fig:proj-filt-Euler-20k}
\end{figure}

At first, we consider a smooth initial condition
\begin{equation}
\label{eq:proj-filt-Euler-exp}
  u_0(x) = \exp(- 20 x^2)
\end{equation}
and simulate in the time interval $[0,4]$ using \num{20000} time steps of the explicit Euler method,
the explicit Euler method with adaptive modal filtering, and the explicit Euler method with a simple
projection.
The simple projection is given by a scaling of all the non-constant Legendre modes by the
same factor, resulting in the desired norm inequality and conservation.
In \autoref{fig:proj-filt-Euler-exp}, we realise that the projection is not
really necessary, the results are very similar to the ones of the
filtered method and all solutions are visually nearly indistinguishable.
Using high order Runge-Kutta schemes does not lead to other observations for
this test case.

However, for the non-smooth initial data
\begin{equation}
\label{eq:proj-filt-Euler-jump}
  u_0(x) =
  \begin{cases}
    1, & -\frac{1}{4} \leq x \leq \frac{1}{4},
    \\
    0, & \text{otherwise},
  \end{cases}
\end{equation}
the same simulation results in Gibbs oscillations and the projection as well as the modal filter have to be applied a
lot more.
The simple projection has also to delete $\Delta t^2 ||\partial_t u||^2$ in every element.
In order to do so, we scale
$\vec{u}_+ =\sum\limits_{n=0}^p u_{+,n} \phi_n$
to
$\vec{u}_+ = u_{+,0} \phi_0 + \alpha
\sum\limits_{n=1}^p u_{+,n} \phi_n$,
where
\begin{equation}
 \alpha
 :=
 \sqrt{
  \frac{ \norm{\vec{u}_+ - u_{+,0} \phi_0}^2 - \Delta t^2 \norm{\partial_t \vec{u}}^2}
        { \norm{\vec{u}_+ - u_{+,0} \phi_0}^2 }
  }
\end{equation}
if $\norm{u_+ - \hat{u}_0\phi_0}^2 - \Delta t^2 \norm{\partial_t u}^2 \geq 0$.
It is not allowed to scale $\hat{u}_0$, since conservation would get lost.
The results of the Euler method using this simple projection fulfilling the
constraints are fairly useless, as can be seen in \autoref{fig:proj-filt-Euler-jump}.
It may be conjectured that the boundary values between cells are influenced in
such a way that the numerical upwind flux adds further dissipation.

\subsubsection*{Simulation using Artificial Viscosity and Modal Filtering}

To validate our investigation from before and especially the adaptive technique and estimation, we
consider the nonlinear Burgers' equation \eqref{eq:Burgers} with smooth initial condition,
\begin{equation}
  \partial_t u + \partial_x \frac{u^2}{2} = 0
  ,\quad
  u(0,x) = u_0(x) = \sin \pi x + 0.01,
\end{equation}
in the periodic domain $x \in [0, 2]$.
This problem serves as a prototypical example of a nonlinear
conservation law, yielding a discontinuous solution in finite time $t \in [0,3]$.
The stable semidiscretisation \eqref{eq:Burgers-semidisc} with $N = 16$ elements
and polynomials of degree $\leq p = 15$
represented w.r.t. a nodal Gauß-Legendre basis is used with the local Lax-Friedrichs flux 
$\fnum(u_-,u_+) = \frac{u_-^2 + u_+^2}{4}
- \frac{ \max \set{ \abs{u_-}, \abs{u_+} } }{2} (u_+ - u_-)$. The explicit Euler
method as time integrator uses $15 \cdot 10^3$ steps for the interval $[0, 3]$.

\begin{figure}[!htb]
  \centering
  ~
  \begin{subfigure}[t]{0.45\textwidth}
    \includegraphics[width=\textwidth]{%
      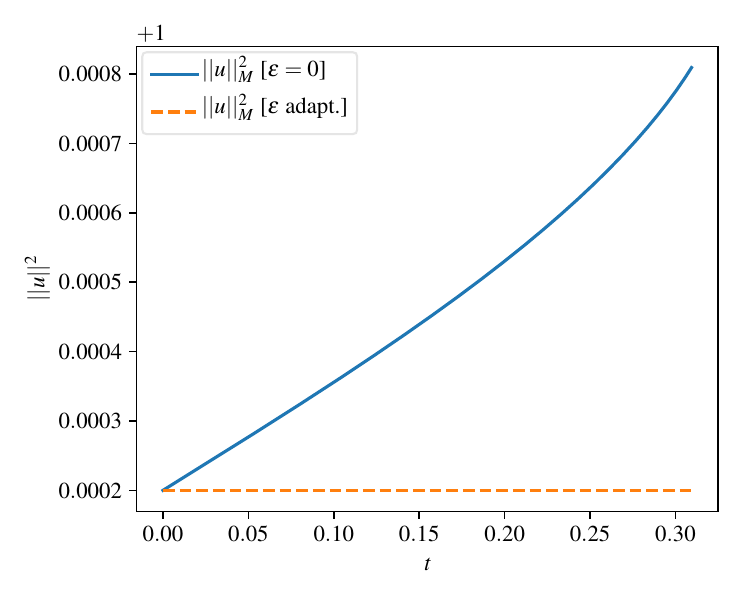}
    \caption{Energy for $t \in [0, 0.31]$. Artificial dissipation.}
    \label{fig:burgers_sin_smooth_energy_av}
  \end{subfigure}%
    ~
  \begin{subfigure}[t]{0.45\textwidth}
    \includegraphics[width=\textwidth]{%
       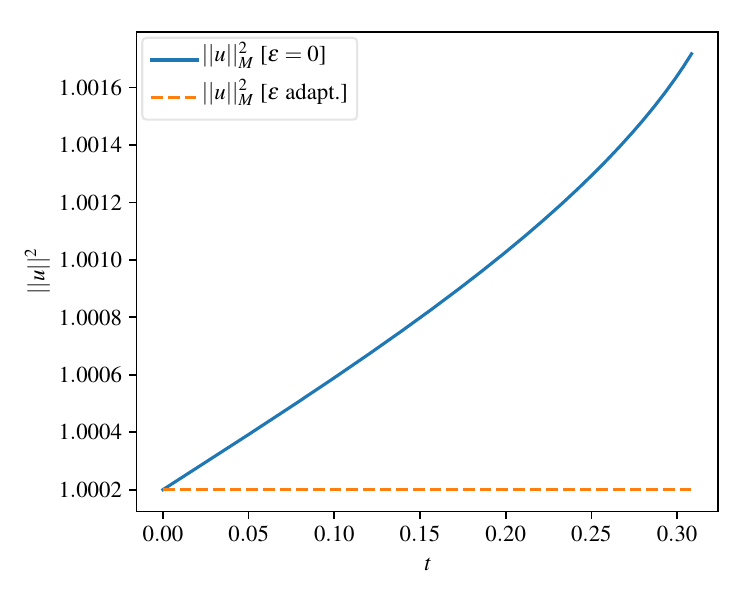}
    \caption{Energy for $t \in [0, 0.31]$. Modal filtering.}
    \label{fig:burgers_sin_smooth_energy_modal}
  \end{subfigure}%
  \\
  ~
  \begin{subfigure}[t]{0.45\textwidth}
    \includegraphics[width=\textwidth]{%
     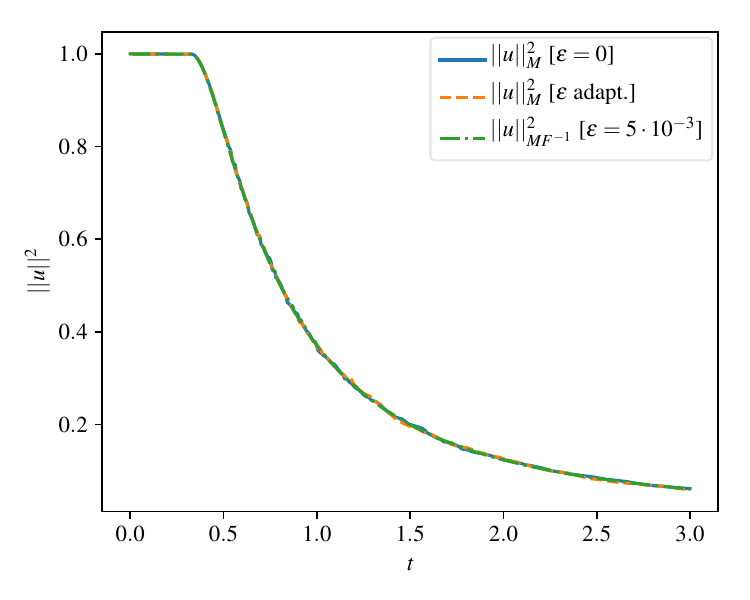}
    \caption{Energy for $t \in [0, 3]$. Artificial dissipation.}
    \label{fig:burgers_sin_rough_energy_av}
  \end{subfigure}%
    ~
  \begin{subfigure}[t]{0.45\textwidth}
    \includegraphics[width=\textwidth]{%
      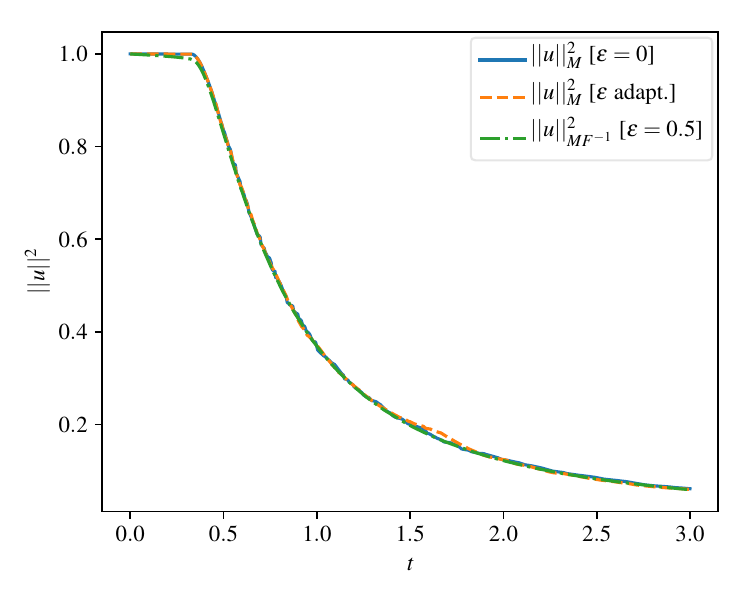}
    \caption{Energy for $t \in [0, 3]$. Modal filtering.}
    \label{fig:burgers_sin_rough_energy_modal}
  \end{subfigure}%
  \caption{Numerical results for Burgers' equation using $N = 16$ elements with
           polynomials of degree $\leq p = 15$.
           The energy are plotted on the left hand side using artificial dissipation and on the right hand side with modal filters.
           }
  \label{fig:burgers_sin}
\end{figure}
First of all, we note that the energy profiles for artificial dissipation (left, i.e. Figure
\ref{fig:burgers_sin_smooth_energy_av} and Figure \ref{fig:burgers_sin_rough_energy_av})
and for modal filtering (right, i.e. Figure
\ref{fig:burgers_sin_smooth_energy_modal} and Figure \ref{fig:burgers_sin_rough_energy_modal})
seem indistinguishable the same.
This demonstrates again that modal filtering can be
seen as the usage of artificial viscosity and vice versa, especially if a similar adaptive strategy is used.

At time $t=0.31$, the solution is still smooth.
However, the energy in Figure \ref{fig:burgers_sin_smooth_energy_av} and
Figure \ref{fig:burgers_sin_smooth_energy_modal} increases if no artificial dissipation or modal filter is applied.
Contrary, applying adaptive artificial dissipation or modal filtering results in a constant energy.
At time $t = 3$, the solution  has developed a discontinuity.
All three choices of artificial dissipation or modal filtering (we compare no filtering, adaptive filtering, and
constant filtering with a fixed strength) yield nearly
visually indistinguishable results for the energy profiles in Figures
\ref{fig:burgers_sin_rough_energy_av} and \ref{fig:burgers_sin_rough_energy_modal} due to the dissipative numerical
flux.

Finally, we would like to mention that around the discontinuities we still obtain oscillations if the adaptive 
strategy is applied
since the strategy is less dissipative and does not smooth out the oscillations from the semidiscrete setting. 
As known in the literature, we can cancel out the oscillations for instance by the usage of limiters 
 which are in accordance with the energy (entropy) inequality \cite{chen2017entropy} or just at the final time by some 
post-processing method.
It should also be noted that the adaptive use of artificial viscosity and modal filtering presented here could 
be used to render shock capturing methods (e.g. \cite{scarnati2018using,glaubitz2019high,glaubitz2019shock, offner2013detecting}) energy 
dissipative, which themselves are not but might have some other advantages.

%% file: 4_explicit_implicit.tex
\section{Comparison Between an Explicit Method with Modal Filtering and the Application of an Implicit Method}
\label{sec:exp_imp}

Here, we present the main part of this work.
It was describe, for instance, in \cite{lozano2018entropy} and \cite{lozano2018entropyImplicit} that
explicit time integration methods may produce entropy whereas in implicit methods entropy may be destroyed.
This entropy production of explicit methods is always a problem when going over from semidiscrete stability to fully
discrete stability.
A classical approach is the usage of implicit methods, for example SBP methods in time, which can be written as
implicit Runge-Kutta methods \cite{nordstrom2013summation, nordstrom2017roadmap, boom2015high}.
Then, the semidiscrete analysis translates directly to the fully discrete scheme.
Unfortunately, this is not the case for explicit methods and in the literature a lot of works can be found which
investigate this issue.
Here, we demonstrate that with our adaptive technique from section \ref{sec:artificial-dissipation}, we can mimic
implicit schemes by using explicit ones with additional dissipation.
As time integration methods we will focus on Euler methods (explicit and implicit).
Further, we will only consider modal filtering, since we have the close connection between modal filtering and
artificial dissipation.

The explicit Euler method
\begin{equation}
  \vec{u}_+ := \vec{u}_0 + \Delta t \, \partial_t \vec{u}_0
\end{equation}
introduces an erroneous growth of energy of size $(\Delta t)^2 \norm{\partial_t \vec{u}_0}^2$,
whereas the implicit Euler method
\begin{equation}
  \vec{u}_+ := \vec{u}_0 + \Delta t \, \partial_t \vec{u}_+
\end{equation}
yields artificial dissipation of size $(\Delta t)^2 \norm{\partial_t \vec{u}_+}^2$
per time step.
Analogously to \ref{subsec:filtering}, the estimate
of the semidiscretisation can be mimicked by filtering with strength
\begin{equation}
  \epsilon
  =
  \left(
    (\Delta t)^2 \norm{ \partial_t \vec{u}_0 }_M^2
  \right)
  \left(
    \sum_{n=0}^p 2 \lambda_n^s \, \tilde{u}_{+,n}^2 \, \norm{\phi_n}^2
  \right)^{-1}
\end{equation}
after each time step.
Similarly, application of this filter and an additional filter with strength
\begin{equation}
  \epsilon
  =
  \left(
    (\Delta t)^2 \norm{ \partial_t \vec{u}_+ }_M^2
  \right)
  \left(
    \sum_{n=0}^p 2 \lambda_n^s \, \tilde{u}_{+,n}^2 \, \norm{\phi_n}^2
  \right)^{-1}
\end{equation}
afterwards yields a filtered explicit Euler method which mimics the dissipation
introduced by an implicit Euler method.
These estimates are applied to the linear advection equation \eqref{eq:transport}
in $[-1,1]$ with periodic boundary conditions.
The initial condition \eqref{eq:proj-filt-Euler-jump} is evolved during the time interval $[0,4]$ on a grid of $N = 8$
elements using polynomials of degree $\leq p = 9$ and an upwind numerical flux.

\begin{figure}[!htb]
  \centering
  \begin{subfigure}[t]{0.45\textwidth}
    \includegraphics[width=\textwidth]{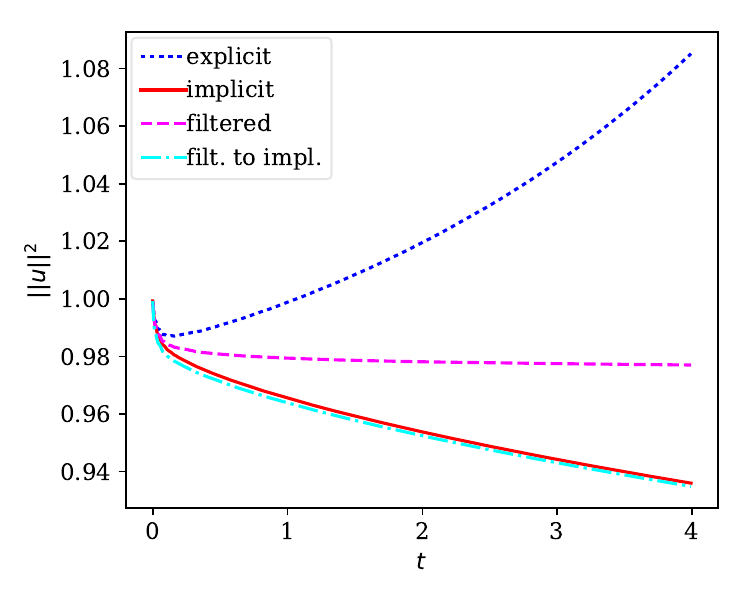}
    \caption{Energy profiles for $t \in [0,4]$ using \num{20000} steps.}
  \end{subfigure}%
    ~
    \begin{subfigure}[t]{0.45\textwidth}
    \includegraphics[width=\textwidth]{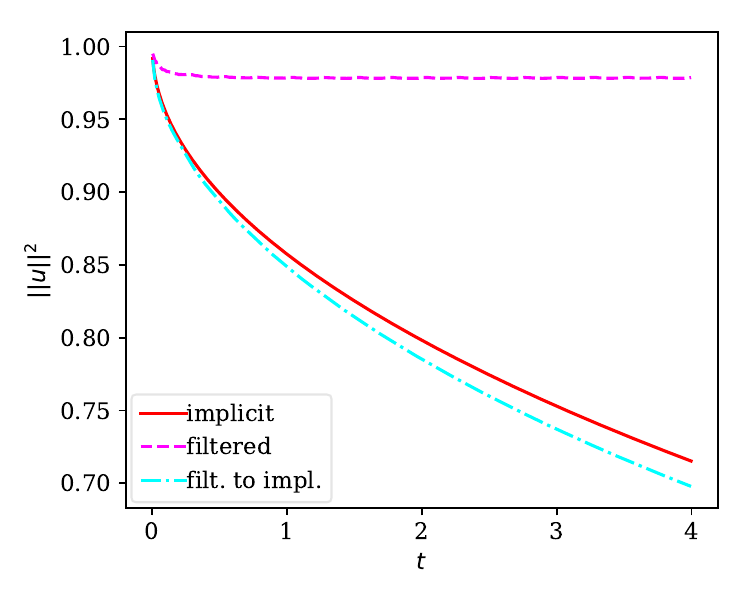}
    \caption{Energies for $t \in [0,4]$ using \num{1000} steps.}
  \end{subfigure}%
  \caption{Energies computed using \num{1000} and \num{20000} time steps of the
           implicit and explicit Euler methods and modal filtering.}
  \label{fig:Euler-20k}
\end{figure}
The corresponding energy profiles using \num{20000} time steps
are plotted in \autoref{fig:Euler-20k}(A) at $t=4$. 
The initial condition \eqref{eq:proj-filt-Euler-jump} is also the exact solution
of the PDE at $t = 4$, i.e. after two periods.
The explicit Euler method (dotted, blue) yields as expected an
unconditional energy growths whereas 
applying adaptive modal filtering once after each time step yields
a nearly constant energy. 
The implicit Euler method (solid red) reduces the energy (introduces artificial
dissipation) as can be seen in the figure.
However, the estimate of the dissipation introduced by implicit Euler yields an  energy result of the explicit
Euler method with modal filtering applied twice (dash-dotted, cyan) that is nearly
indistinguishable from the implicit one.

%

Although the estimate of the filter strength is conservative (i.e. only necessary),
the energy of the twice filtered explicit solution is slightly less than the energy
of the implicitly computed solution. The reason is probably the appearance of some changes of
boundary values due to the filtering that triggers additional dissipation by the
upwind flux.


Finally, we note that the same behaviour can be observed  if one uses considerably
less time steps.
In \autoref{fig:Euler-20k}(B) the results of the implicit and filtered explicit Euler method using only
\num{1000} time steps are plotted.
Similarly to the case before, the filtered solutions approximate their targets very well.

%% file: 5_Summary.tex
\section{Summary}
\label{Summary}

The application of SBP operators in time together with a semidiscrete method
yields to a fully discrete stable scheme to solve hyperbolic conservation laws as it is done in
\cite{nikkar2015fully, nordstrom2017roadmap, friedrich2018entropy} 
whereas in \cite{ketcheson2019relaxation, ranocha2019relaxation} a relaxation RK method is applied.
Here, we follow another approach and consider also a semidiscretely stable scheme and \textbf{explicit} time integration methods
but to reach a fully discrete stable scheme, we apply modal filtering or artificial dissipation,
where the strength of dissipation is steered automatically by an adaptive strategy.
We consider only the explicit Euler method in this context.
However, since strong stability preserving Runge-Kutta schemes can be written as convex combination of explicit
Euler steps, our approach can be extended to these methods.
Then, we demonstrated by a concrete example that with the usage of modal filters together with our adaptive strategy,
we are able to mimic the behavior of an implicit method and can imitate the stability properties of this scheme.
This contribution leads to a better understanding of existing algorithms and numerical techniques, especially the
application of artificial dissipation as well as modal filtering in the context of numerical methods for hyperbolic
conservation laws together with the selection of explicit or implicit time integration methods.
A future research topic will be further extension of the study presented here together with implicit SBP operators
in time.
Also the usage of other adaptive strategies, such as the annihilation of the entropy production in time
\cite{lozano2018entropy, lozano2018entropyImplicit} will be considered.

%% file: Appendix.tex
\section*{Appendix}
\label{sec:appendix}

In this section, we show how our analysis from subsection \ref{subsec:filtering}
can be applied to RK methods and transfer the results to DEC methods, which can be formulated in the RK framework 
as well. 
A RK method with $s$ stages is  given by its
Butcher tableau
\begin{equation}
\label{eq:butcher}
\begin{array}{c | c}
  c & A
  \\ \hline
    & b
\end{array}.
\end{equation}
Here, $A \in \R^{s \times s}$ and $b, c \in \R^s$.
Since there is no explicit dependence on the time in the semidiscretisation,
one step from $\vec{u}_0$ to $\vec{u}_+$ is given by
\begin{equation}
\label{eq:erk}
  \vec{u}_i
  :=
  \vec{u}_0 + \Delta t \sum_{j=1}^{s} a_{ij} \, \partial_t \vec{u}_j,
  \qquad
  \vec{u}_+
  :=
  \vec{u}_0 + \Delta t \sum_{i=1}^{s} b_{i} \, \partial_t \vec{u}_i.
\end{equation}
Here, the $\vec{u}_i$ are the stage values of the RK method. 
It is also possible to express the method via the slopes $\vec{k}_i = \partial_t \vec{u}_i$,
as done  by \cite[Definition II.1.1]{hairer2006geometric}. 
Using the stage values $\vec{u}_i$ as in \eqref{eq:erk}, we have 
\begin{align*}
\stepcounter{equation}\tag{\theequation}
\label{eq:estimate-RK}
  &
  \norm{ \vec{u}_+ }_M^2 - \norm{ \vec{u}_0 }_M^2
  \\=&
  2 \Delta t \scp{ \vec{u}_0 }{ \sum_{i=1}^{s} b_{i} \, \partial_t \vec{u}_i }_M
  + (\Delta t)^2 \norm{ \sum_{i=1}^{s} b_{i} \, \partial_t \vec{u}_i }_M^2
  \\\stackrel{\eqref{eq:erk}}{=}&
  2 \Delta t \sum_{i=1}^{s} b_{i}
    \scp{ \vec{u}_i - \Delta t \sum_{j=1}^{s} a_{ij} \, \partial_t \vec{u}_j }{
    \partial_t \vec{u}_i }_M
  + (\Delta t)^2 \norm{ \sum_{i=1}^{s} b_{i} \, \partial_t \vec{u}_i }_M^2
  \\=&
  2 \Delta t \sum_{i=1}^{s} b_{i} \scp{ \vec{u}_i }{ \partial_t \vec{u}_i }_M
  + (\Delta t)^2 \left[
    \norm{ \sum_{i=1}^{s} b_{i} \, \partial_t \vec{u}_i }_M^2
    - 2 \sum_{i,j=1}^{s} b_{i} \, a_{ij}
      \scp{ \partial_t \vec{u}_i }{ \partial_t \vec{u}_j }_M
  \right]
  \\=&
  2 \Delta t \sum_{i=1}^{s} b_{i} \scp{ \vec{u}_i }{ \partial_t \vec{u}_i }_M
  + (\Delta t)^2 \left[
    \sum_{i,j=1}^{s} \left( b_i b_j - b_{i} \, a_{ij} - b_j a_{ji} \right)
      \scp{ \partial_t \vec{u}_i }{ \partial_t \vec{u}_j }_M
  \right],
\end{align*}
where the symmetry of the scalar product has been used in the last step.
Here, the first term on the right hand side is consistent with
$\int_{t_0}^{t_0 + \Delta t} 2 \scp{u}{\partial_t u}$, if the RK method
is consistent, i.e. $\sum_{i=1}^{s} b_i = 1$. 

The second term of order $(\Delta t)^2$ is undesired. 
Depending on the method (and the stages, of course), it may be positive or negative. 
However, if it is positive, then a stability error may be introduced. 
As a special case, if the method fulfils $b_i b_j = b_i a_{ij} + b_j a_{ji},\, i,j \in \set{1, \dots, s}$,
this term vanishes. 
Such methods can conserve quadratic invariants of ordinary
differential equations, a topic of geometric numerical integration, see
Theorem IV.2.2 of \cite{hairer2006geometric}, originally proven by
\cite{cooper1987stability}. 
A special kind of these methods are the implicit Gauss methods, see section II.1.3 of \cite{hairer2006geometric}.

For an explicit method ($a_{ij} = 0$ for $j \geq i$), the undesired term of order
$(\Delta t)^2$ in \eqref{eq:estimate-RK} can be rewritten as
\begin{equation}
\label{eq:Delta-t-square-term-explicit}
\begin{aligned}
  &
  \norm{ \sum_{i=1}^{s} b_{i} \, \partial_t \vec{u}_i }_M^2
  - 2 \sum_{i=1}^{s} \sum_{j=1}^{i-1}  b_{i} \, a_{ij}
    \scp{ \partial_t \vec{u}_i }{ \partial_t \vec{u}_j }_M
  \\=&
  \sum_{i=1}^{s} b_{i}^2 \norm{\partial_t \vec{u}_i}_M^2
  + 2 \sum_{i=1}^{s} \sum_{j=1}^{i-1}  b_{i} \, (b_j - a_{ij})
    \scp{ \partial_t \vec{u}_i }{ \partial_t \vec{u}_j }_M.
\end{aligned}
\end{equation}

This undesired increase of the norm may be remedied by the application of
an adaptive modal filter $\mat{F}$. Analogously to subsection \eqref{subsec:filtering}, 
the adaptive filter strength $\epsilon$ may be estimated via
\begin{equation}
\begin{aligned}
  \norm{ \mat{F} \vec{u}_+ }_M^2
  &\stackrel{!}{\leq}
  RHS
  :=
  \norm{ \vec{u}_0 }_M^2
  + 2 \Delta t \sum_{i=1}^{s} b_{i} \scp{ \vec{u}_i }{ \partial_t \vec{u}_i }_M
  \\&\leq
  RHS
  + (\Delta t)^2 \left[
    \sum_{i,j=1}^{s} \left( b_i b_j - b_{i} \, a_{ij} - b_j a_{ji} \right)
      \scp{ \partial_t \vec{u}_i }{ \partial_t \vec{u}_j }_M
  \right],
\end{aligned}
\end{equation}
if the term of order $(\Delta t)^2$ is non-negative.
In a modal Legendre basis $\set{\phi_i}$, the (exact) modal filter $\mat{F}$ is given by \eqref{eq:modal_filter}. 
Thus,
\begin{equation}
  \sum_{n=0}^p \exp[- 2 \epsilon \, \lambda_n^s] \, u_{+,n}^2 \, \norm{\phi_n}^2
  \stackrel{!}{\leq}
  RHS 
\end{equation}
is required. 
Here, $u_{+,n}$ are the coefficients of the polynomial $u_+$, expressed in the
Legendre basis of polynomials of degree $\leq p$. 
Following \eqref{eq:approximation}, the filter strength $\epsilon$ can be estimated by
\begin{equation}
\begin{aligned}
  &
  \sum_{n=0}^p (1 - 2 \epsilon \, \lambda_n^s) \, u_{+,n}^2 \, \norm{\phi_n}^2
  \leq RHS
  \\ \iff
  &
  \left( \sum_{n=0}^p u_{+,n}^2 \, \norm{\phi_n}^2 - RHS \right)
  \left( \sum_{n=0}^p 2 \lambda_n^s \, u_{+,n}^2 \, \norm{\phi_n}^2 \right)^{-1}
  \leq
  \epsilon,
\end{aligned}
\end{equation}
for $\sum_{n=0}^p 2 \lambda_n^s \, u_{+,n}^2 \, \norm{\phi_n}^2 > 0$.
Using $\sum_{n=0}^p u_{+,n}^2 \, \norm{\phi_n}^2 \approx \norm{ \vec{u}_+ }_M^2$
(since $\norm{\cdot}_M$ approximates the exact $\L_2$ norm on the left hand side),
this results in
\begin{equation}
\label{eq:estimate-epsilon_RK}
  \epsilon
  \geq
  \left(
    \norm{ \vec{u}_+ }_M^2 - \norm{ \vec{u}_0 }_M^2
    - 2 \Delta t \sum_{i=1}^{s} b_{i} \scp{ \vec{u}_i }{ \partial_t \vec{u}_i }_M
  \right)
  \left(
    \sum_{n=0}^p 2 \lambda_n^s \, \tilde{u}_{+,n}^2 \, \norm{\phi_n}^2
  \right)^{-1}.
\end{equation}
\eqref{eq:estimate-epsilon_RK}  is the general estimation. 
We can derive estimation \eqref{eq:estimate-epsilon} from \eqref{eq:estimate-epsilon_RK} using the coefficients for 
the explicit Euler method.